\documentclass[UTF8,10pt,fontset=windows,twoside]{article} %若是windows系统下TexLive2018及以前的版本，则可以去掉fontset=windows这一段命令。
\usepackage[left=3.18cm,right=3.18cm,top=2.54cm,bottom=2.54cm]{geometry}
\usepackage{amsfonts}
\usepackage{amsmath}
\usepackage{amssymb}
\usepackage{mathtools}
\usepackage{setspace}
\usepackage{amsthm}
\usepackage{pgfplots}
\usepackage{enumitem}
\usepackage{tikz}   %交换图表
\usepackage{tikz-cd}  %交换图表
\usepackage{changepage}
\usetikzlibrary{arrows}
\usepackage{subcaption}
\usepackage{titlesec}
\usepackage{floatrow}
\usepackage{bm}
\usepackage{graphicx}
\usepackage{cases}
\usepackage{appendix}
\usepackage[colorlinks,linkcolor=black,anchorcolor=black,citecolor=black]{hyperref}
\usepackage[numbers]{natbib}
\usepackage{txfonts}  %请将该数学字体宏包置于各种数学公式宏包后面，否则容易产生冲突
\usepackage{fancyhdr}
\fancypagestyle{plain}{}
\title{A Note Related to Graph Theory}
\author{Jinfeng Li}
\usepackage[scaled=0.92]{helvet}	% set Helvetica as the sans-serif font
\theoremstyle{definition}
\newtheorem{thm}{theorem}[section]

\newtheorem{claim}[thm]{claim}
\newtheorem*{thmm}{theorem}  %不带编号
\newtheorem{obs}{observation}

\newtheorem{lem}[thm]{lemma}

\newtheorem{conj}{conjecture}[section]
\newcommand{\pf}{proof}  
\numberwithin{equation}{section}
\def\pf{\noindent {\it Proof. }}
\usepackage{xcolor}
\usepackage{listings}

\begin{document}

    \maketitle
\begin{abstract}
        This article foucuses on $(P_3\cup P_2,K_4)$-free graph. In this paper, we prove that if G is $(P_3\cup P_2,K_4)$-free, then $\chi(G)\le 7$. We then use this result to obtain the upper bound of order and chromatic number of $(4K_1,\overline{P_3\cup P_2},K_{\omega})$-free graph .
\end{abstract}

\section{Introduction}

\quad A graph $G$ consists of its vertex set $V(G)$ and edge set $E(G) \subseteq V(G)\times V(G)$. The {\bf order} of $G$, denoted by $n$, is the size of $V(G)$.  The {\bf complement} of a graph $G$, denoted by $\overline{G}$ or $co-G$, is a graph with the same vertex set while whose edge set consists of the edges not present in $E(G)$.

   All graphs in this paper are finite and have no loops or parallel edges. A graph $G$ {{\bf contains} $H$ if $H$ is isomorphic to an induced subgraph of $G$. If $H_1$ is {\bf isomorphic} to $H_2$, then we say $H_1\simeq H_2$. For a collection of graphs $H_t,t=1,2,...,n$, $G$ is $(H_1,H_2,...,H_n)$-free if it does not contain $H_t,t=1,2,...,n$. In this paper, if $G$ does not contain $P_2$, we will say that $G$ is {\bf edge-free}.

    Path and cycle on $n$ vertex is denoted by $P_n$ and $C_n$, respectively. The complete graph on that has $n$ vertex is denoted bu $K_n$. We use $G\cup H$ to denote the disjoint union of $G$ and $H$.   

   A {\bf $k$-coloring} of a graph $G$ is a function $\phi:V(G)\rightarrow \{1,...,k\}$ such that $\phi(u)\ne \phi(v)$ whenever $u$ and $v$ are adjacent in $G$. We say that $G$ is $k$-colorable if $G$ admits a $k$-coloring. The {\bf chromatic number} of $G$ is denoted by $\chi(G)$ which represents the minimum positive integer $k$ such that $G$ is $k$-colorable. The {\bf clique number} is denoted by $\omega(G)$ which represents the size of the largest clique in $G$.

  A family $\mathbb{G}$ of graphs is said to be {\bf $\chi$-bounded} if there exists a function f such that for every graph $G \in \mathbb{G}$ and every induced subgraph $H$ of $G$ it holds that $\chi(H)\le f(\omega(H))$. The function $f$ is called a $\chi$-binding function for $G$. The class of perfect graphs (a graph $G$ is perfect if for every induced subgraph $H$ of $G$ it holds that $\chi(H) = \omega(H))$, for instance, is a $\chi$-bounded family with $\chi$-binding function $f(x) = x$. Therefore, $\chi$-boundedness is a generalization of perfection. The notion of $\chi$-bounded families was introduced by Gy\'arf\'as who make the following conjecture.

\begin{conj}{\cite{G1973}}
  For every forest T, the class of $T$-free graphs is $\chi$-bounded.
\end{conj}

  For $\omega=3$, Esperet, Lemoine, Maffray and Morel{\cite{EL2013}} obtained the optimal bound on the chromatic number: every $(P_5, K_4)$-free graph is 5-colorable. Serge Gaspers and Shenwei Huang{\cite{SH2018}} obtained the optimal bound of the chromatic number: every $(2P_2,K_4)$-free graph is 4-colorable. Karthick and S. Mishra obtained the optimal bound of the chromatic number: every $(P_6,diamond,K_4)$-free graph is 6-colorable. Rui Li,Jinfeng Li and Di Wu claimed that {\cite{LD2023}} $(P_3\cup P_2,K_4)$-free graph is $9$-colorable. We follow the methods in this paper and sharper the bound to $7$.  

   Partition $V(G)$ into two following parts:

\begin{itemize}
   \item $D_1:=\{x\in V(G)|\omega(G-N(x))\le 2  \}$
   \item $D_2:=G-D_1$ 
\end{itemize}

   Let $C = u_1v_1v_2v_3u_3u_2$ be a $6$-hole. We use {\bf $co-donino$} to denote a graph obtained from $C$ by connecting edges $v1v3$ and $u1u3$.

   Let $C = u_1v_1v_2v_3u_3u_2$ be a $6$-hole. We use {\bf $co-A$} to denote a graph obtained from $C$ by connecting edges $v_1v_3$  ,$u_1u_3$ and $u_1v_3$.

   Our proof includes three major theorem:

\begin{thm}{\label{t1}}
   If $G[D_1]$ contains $P_2\cup P_1$, then $\chi(G)\le 7$.
\end{thm}

\begin{thm}{\label{t2}}
   If $G$ contains $co-domino$ or $co-A$, then $\chi(G)\le 7$.
\end{thm}

\begin{thm}{\label{t3}}
   If $G$ is $(co-domino,co-A)$-free and $G[D_2]$  is not an induced subgraph of $K_3$, then $\chi(G)\le 7$. 
\end{thm}

   Finally, we obtain a theorem for $(4K_1,\overline{P_3\cup P_2})$-free graphs:
\begin{thm}
   If $G$ is $(4K_1,\overline{P_3\cup P_2})$-free with clique number $\omega$ and order n, then $n\le 7\omega$ and $\chi(G)\le 4\omega$.
\end{thm}

\section{Structure Around $K_3$}    
   
\quad       Suppose there is a $K_3$ in $G$ and its vertex set is $\{v_1,v_2,v_3\}$. Naturally, we can divide $V(G)$
into following three sets:
  \begin{itemize}
\item $A_0:=\{v|v\text{ is not adjacent to any vertex in } \{v_1,v_2,v_3\}\}$.
\item $A_1:=\{v|v\text{ is adjacent to one vertex in } \{v_1,v_2,v_3\} \text{ only}\}$.
\item $A_2:=\{v|v\text{ is adjacent to two vertex in } \{v_1,v_2,v_3\} \text{ only}\}$.
\end{itemize}

  we need to combine $A_0$ and $A_1$ and divide it into three disjoint sets: $B_1,B_2,B_3$. The set $B_1$ includes vertices that is only adjacent to $u_1$ and vertices in $A_0$. The set $B_j,j=2,3$ includes vertices that are only adjacent to $u_j$. Since $G[B_i]$ is anticomplete to an edge in $G[\{v_1,v_2,v_3\}]$, $G[B_i],i=1,2,3$ is $P_3$-free. 

  $A_2\cup \{v_1,v_2,v_3\}$ can be partitioned in to  $(N(v_i)\cap N(_{i-1}))\cup \{v_{i+1}\}$(mod3) for $i=1,2,3$. It is easy to obtain that $\chi((N(v_i)\cap N(_{i-1}))\cup \{v_{i+1}\})\le 1$(mod3) for $i=1,2,3$.

\section{main theorem}

\quad We first introduce a lemma provided by Rui Li, Jinfeng Li and Di Wu.

\begin{lem}{\cite{LD2023}}
  If $G$ contains $P_2\cup K_3$, then $\chi(G)\le 6$.
\end{lem}

       Define $D_1:=\{x|\omega(G-N(x))\le 2\}$ and $D_2:=G-D_1$.

\begin{claim}{\label{c1}}
    For any two nonadjacent vertex  $y_1,y_2$ in $D_1$, $\omega(G[N(y_1)-N(y_2)])\le 1$ and $\omega(G[N(y_2)-N(y_1)])\le 1$. 
\end{claim}       
\pf
        If there is an edge in $G[N(y_1)-N(y_2)]$(or $G[N(y_2)-N(y_1)]$), then $\omega(G[N(y_1)-N(y_2)])=3$(or $\omega(G[N(y_2)-N(y_1)])=3$), which contradicts definition of $y_1$(or $y_2$).\qed

       We fisrt introduce a lemma to help us slightly determine the structure of $G[D_1]$. 

\begin{thmm}{(1.1)}
      If $G[D_1]$ contains $P_2\cup P_1$, then $\chi(G)\le 7$.
\end{thmm}
\pf
     Suppose $G[\{v_1,v_2\}\cup \{v_3\}]\simeq P_2\cup P_1$, $M(v_1,v_2)=G-\{v_1,v_2\}-(N(v_1)\cup N(v_2))$. we can divide $G$ into $\{v_1,v_2\}\cup (N(v_1)\cup N(v_2))\cup M(v_1,v_2)$. Since $G$ is $(P_3\cup P_2,P_2\cup K_3)$-free, $G[M(v_1,v_2)]$ is unions of clique and $\omega(M(v_1,v_2))\le 2$. Therefore, $\{v_1,v_2\}\cup M(v_1,v_2)\le 2$.

     We will prove $\chi(N(v_1)\cup N(v_2))\le 5$ and hence $\chi(G)\le 7$. $N(v_1)\cup N(v_2)$ can be divided into $N(v_1)-N(v_3)$, $(N(v_1)\cap N(v_3))-N(v_2)$, $N(v_2)-N(v_3)$, $(N(v_2)\cap N(v_3))-N(v_1)$ and $N(v_1)\cap N(v_2)$. We apply claim \ref{c1} to prove that  the  clique number of  first four sets are at most 1. If the last set has edge, then $\omega(G)\ge 4$. Therefore, $\chi(N(v_1)\cup N(v_2))\le 5$.  \qed

\begin{obs}
   If $G[D_1]$ is $P_2\cup P_1$-free, then $\overline{G[D_1]}$ is $P_3$-free.
\end{obs} 
 
  Now we introduce a lemma which is important in the proof of theorem1.2.

\begin{claim}{\label{c1}}
    Suppose $V(G)$ can be partitioned  into three nonempty subsets $V_1$ , $V_2$ and $V_3$ such that $G[V_1]$ and $G[V_2] $  are $(K_3,P_3)$-free graph,  $G[V_3]$ is a stable set , $G[V_1 \cup  V_2]$ is  $K_3$-free , $G[V_j \cup V_3]$ is $(K_3,P_3)$-free graph  for $j=1,2$ and any vertex $v$ in $V_i,i=1,2$, $M(v)-V_i$ has no edge. Furthermore, for any vertex $v\in V_3$, if $N(v)\cap V_i\ne \emptyset$, then $N(v)\cap V_i=\{x\in V_i| x\text{ is not complete to } N(v)\cap N_{3-i}\},i=1,2$. If $G[V_1]$ has at most one edge, then $\chi(G)\le 3$.
\end{claim}

\proof
 
  If $V_1$ has no edge, then randomly select one vertex $v\in V_1$. We partition $V(G)$ accoding to $\{v\}$, that is $V(G)=\{v\}\cup N(v) \cup M(v)$.  $|N(v)\cap V_3|\le 1$. Now we partition $V(G)$ into three stable set:  $(N(v)\cap V_3) \cup (V_1\cap M(v))$, $N(v)\cap V_2$ and  ${v}\cup (M(v)-V_1)$. $(N(v)\cap V_3)\cup (V_1\cap M(v))$   is stable set as $G[V_3\cup V_1]$ is $P_3$-free.  According to definition, the rest two sets are stable sets.

 If $V_1$  has one edge, then set one vertex of the vertex set of that edge $v$. We partition $V(G)$ into three part: $V(G)=\{v\}\cup N(v)\cup M(v)$. It is obvious that $N(v)\subseteq V_2$, otherwise $G[V_1\cup V_2]$ induces an $P_3\cup P_2$ or $K_3\cup P_2$. Therefore, as $G[V_1\cup V_2]$ is $K_3$-free, we can conclude that $\chi(N(v))\le 1$.   $\chi(M(v))\le \chi(M(v)\cap V_1)+\chi( M(v)-V_1)\le 2$ and hence $\chi(G)\le \chi(N(v))+\chi(\{v\}\cup M(v))\le 3$. \qed

\begin{obs}
   If $G$ contains $co-domino$, then according to section 2.1, $V(G)$ can be partitioned into $B_1,B_2,B_3$ and $A_2$.
\end{obs}
      
\begin{lem}{ \label{l1}}
    If $G$ contains  $co-domino$, then $\chi(G)\le 7 $.
\end{lem}
     
\pf

 \begin{figure}[H]
\centering
\includegraphics[width=0.2\textwidth]{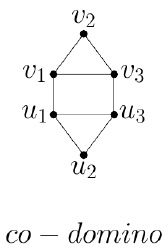}

\end{figure}    

 $\omega(G[B_2])\le 1$. Otherwise there exists an edge in $G[B_2]$. We call the edge $G[\{y_1,y_2\}]$. If $y_1\sim u_2$, then $G[\{v_1,v_3\}\cup \{y_2,y_1,u_2\}]\simeq P_3\cup P_2$~or~$G[\{v_1,v_3\}\cup \{y_2,y_1,u_2\}]\simeq K_3\cup P_2$(depending on whether $y_2$ is adjacent to $u_2$). As a consequence, both $y_1$ and $y_2$ are not adjacent to $u_2$. However, $y_1$ is complete to $\{u_1,u_3\}$, which can be proven through making use of $P_2\cup P_3$-free. Similarly, $y_2$ is complete to $\{u_1,u_3\}$. As a consequence, $G[\{y_1,y_2,u_3,u_1\}]\simeq K_4$, which contradicts the fact that $\omega(G)\le 3$. As a consequence, $\chi(B_2)\le 1$.

   For any $x\in B_2$, $N(x)$ must be one of four following sets:  $\{u_1,u_2\}$,$\{u_2,u_3\}$,$\{u_2\}$ and $\{u_3,u_1\}$. We partition $B_2$($\{u_1,u_2\}$) into $C_1$($\{u_2,u_3\}$),$C_2$($\{u_2,u_3\}$),$C_3$($\{u_2\}$) and $C_4$($\{u_3,u_1\}$) accoding to their neigbors in $\{u_1,u_2,u_3\}$. Trivially  , $C_i\cap C_j=\empty$ for $i\ne j\in \{1,2,3,4\}$.

    If $C_1\ne \emptyset$, then  there is no edge in $G[B_3]$.  Suppose $\{y_1,y_2\}\in E(B_3)$ and $x\in C_1$, we are going to prove $y_i,i=1,2$ are complete to $\{v_1,x\}$ and anticomplete to $\{u_2,u_3\}$ and hence $G[\{y_1,y_2,u_1,x\}]\simeq K_4$, which contradicts $\omega(G)\textless 4$. If $y_1$ is adjacent to $u_2$, then $G[\{y_1,u_3,u_2\}\cup \{v_2,v_1\}]$ is isomorphic to $P_2\cup P_3$ or $P_2\cup K_3$. Symmetrically, $y_2$ is anticomplete to $\{u_2,u_3\}$. As consequence, $y_1$ must be complete to $\{u_1,x\}$, or $G[\{u_2,u_1,x,v_2,v_3,y_1\}]$ will induce $P_3\cup P_2$. Symmetrically, $y_2$ is complete to $\{x,u_1\}$. As consequence, we get the contradiction that $G[\{y_1,y_2,u_1,x\}]\simeq K_4$. Therefore, $\omega(B_3)\le 1$.

   If $C_2\ne \emptyset$, then there is no edge in $G[B_1-A_0]$, which is similar to the proof for $C_1\ne \emptyset.$ 

   From dicussion above, if $C_1\cup C_2\ne \emptyset$, then $\chi(B_1\cup B_2\cup B_3)\le 4$ and hence $\chi(G)\le \chi(B_1\cup B_2\cup B_3)+\chi(A_2\cup \{v_1,v_2,v_3\})\le 4+3=7$. 

 \begin{spacing}{2.0}
\end{spacing}

   Suppose $C_1\cup C_2=\emptyset$ and $x\in C_3$. Furthermore, $|C_3|\le 1$. Otherwise, $G[\{u_2\}\cup \{v_1,v_3\}\cup C_3]$ induces $P_3\cup P_2$ or $P_2\cup K_3$. According to the definition, $A_2$ can be partitioned  into $N(v_1)\cap N(v_2)$, $N(v_1)\cap N(v_3)$ and $N(v_2)\cap N(v_3)$. 

   Suppose there is a vertex  in $(N(v_2)\cap N(v_3)$ has only one neighbor in $\{u_1,u_2,u_3\}$. It is not difficult to obtain that the only neighbor is $v_2$. If we exchange $\{v_1,v_2,v_3\}$ and $\{u_1,u_2,u_3\}$ in $co-domino$, then we can prove $\chi(G)\le 7$ in the same way as what we did when $C_1\cup C_2\ne \emptyset$.  

   Suppose evrey vertex in $A_2\cap N(v_2)$ has two neighbors in $\{u_1,u_2,u_3\}$. 

   Let $C = u_1v_1v_2v_3u_3u_2$ be a $co-domino$ and $\{u,x\}$ be two vertex outside $C$. We use {\bf co-donino1} to denote a family of graphs obtained from $\{u,x\}\cup C$ by connecting edges $uu_1,uu_2,uv_1,uv_2,xv_2,xu_2$ and $xu$ can be connected or disconnected.

   Let $C = u_1v_1v_2v_3u_3u_2$ be a $co-domino$ and $\{u,x\}$ be two vertex outside $C$. We use {\bf co-domino2} to denote a graph obtained from $C$ by connecting edges $uu_1,uu_2,uv_1,uv_3,xv_2,xu_2$ and $xu$.
 \begin{spacing}{2.0}
\end{spacing}
   Suppose $u \in N(v_2)\cap N(v_1)$.  The following case shows that if $u$ is complete to $\{u_1,u_2\}$ or complete to $\{u_2,u_3,x\}$, then $\chi(G)\le7$. In other words, if $G$ contains co-domino1 or co-domino2, then $\chi(G)\le 7$.  
\begin{figure}[H]
\centering

\includegraphics[width=0.30\textwidth]{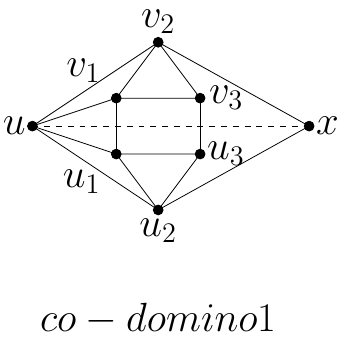}
\includegraphics[width=0.30\textwidth]{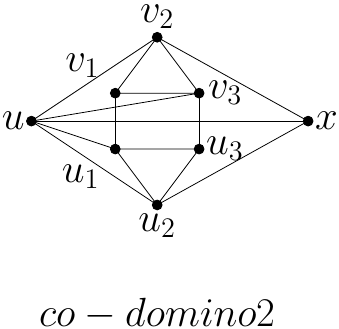}

\end{figure}    
Recalling that $G[B_3]$ is anticomplete to $u_2$ ,  $B_3-\{u_3\}$ is anticomplete to $u_2$.
\begin{spacing}{2.0}

\end{spacing}
{\bf case1:} $G$ contains $co-domino1$ or $co-domino2$.
   
   Suppose $G$ contains $co-domino1$.

    $G[B_3]$ contains at most one  edge. Otherwise, suppose $G[B_3]$ contains more than one edge. Since $G[B_3]$ is $P_3$ -free and complete to $\{u_1\}$ and anticomplete to $\{u_3\}$,  there exists an $m$ such that these edges are isomorphic to $mK_2,m\textgreater 1$ and $u$ has at most one neighbor in every edges or $G[\{u,u_1\}\cup B_3]$ induces $K_4$. However, $m\textgreater 1$ leads to a contradiction that  $G[\{B_3-N(u)\}\cup \{u,u_2\}]$ induces $ P_3\cup P_2$. 
     
    Let  $G[B_1-A_0]=V_1',G[B_3]=V_2',G[A_0]=V_3'$. According to claim\ref{c1}, $\chi(B_3\cup B_1)\le 3$. Therefore, $\chi(G)\le \chi(B_1\cup B_3)+\chi(B_2)+\chi(A_2)\le 3+1+3=7$. 

\begin{spacing}{2.0}
\end{spacing}

    Suppose $G$ is $co-domino1$-free and contains $co-domino2$.

    $G[B_1-A_0]$  contains at most one  edge. Since $G[B_1-A_0]$ is $P_3$ -free and complete to $\{u_1\}$ and anticomplete to $\{u_3\}$,    this union of edges is isomorphic to $mK_2,m\textgreater 1$ for some $m$ and $u$ has at most one neighbor in every edge or $G[\{u,x\}\cup B_3]$ induces $K_4$. However, $m\textgreater 1$ leads to a contradiction that  $G[\{B_3-N(u)\}\cup \{u,u_2\}]$ induces $ P_3\cup P_2$. 

  Let $G[B_1-A_0]=V_1',G[B_3]=V_2',G[A_0]=V_3'$, then according to claim\ref{c1}, $\chi(B_3\cup B_1)\le 3$. Furthermore, $\chi(G)\le \chi(B_1\cup B_3)+\chi(B_2)+\chi(A_2)\le 3+1+3=7$. 
\begin{spacing}{2.0}
\end{spacing}
   Let Let $C = u_1v_1v_2v_3u_3u_2$ be a $co-domino$ and $\{u,x\}$ be two vertex outside $C$. We use {\bf co-domino3} to denote a graph obtained from $C$ by connecting edges $uu_1,uu_2,uv_1,uv_3,xv_2$ and $xu_2$.
\begin{figure}[H]
\centering

\includegraphics[width=0.30\textwidth]{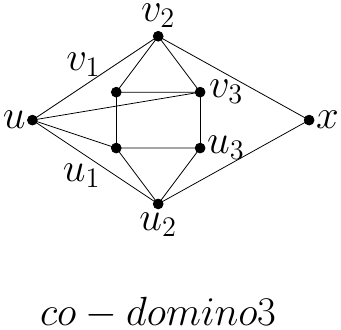}

\end{figure} 
Suppose $u\in N(v_1)\cap N(v_2)$ and $u\sim u_2,u\sim u_3,u\not\sim x$.  The following case shows that if $u$ exists, then $\chi(G)\le 7$. In other words, if $G$ contains $co-domino3$, then $\chi(G)\le 7$.

 Partition $A_2=(N(v_1)\cap N(v_2))\cup(N(v_2)\cap N(v_3))\cup (N(v_1)\cup N(v_3))$ into following five vertex sets.

\begin{centering}
$D_1:=N(v_1)\cap N(v_2)\cap N(u_2)\cap N(u_3)$.\\
$D_2:=N(v_1)\cap N(v_2)\cap N(u_1)\cap N(u_3)$.\\
$D_3:=N(v_2)\cap N(v_3)\cap N(u_1)\cap N(u_2)$.\\
$D_4:=N(v_2)\cap N(v_3)\cap N(u_1)\cap N(u_3).$\\
$D_5:=N(v_1)\cap N(v_3).$\\
\end{centering}

{\bf case2:} Suppose $G$ contains $co-domino3$.

    Suppose $\{u,u'\}\subseteq E[D_1,D_3]$.

    $\chi(B_1\cup B_3-\{u_1,u_3\})\le 3$.  Define $I:=B_1\cup B_3-\{u_1,u_2,u_3\}$. Recalling that $B_1\cup B_3-A_0$ is anticomplete to $\{u_2\}$ and that $A_0$ is anticomplete to $\{v_1,v_2,v_3\}$, $I$ is anticomplete to $\{v_2,u_2\}$. Therefore, $I-N(u)$ is anticomplete to $\{v_2,u,u_2\}$. Since $G[\{v_2,u,u_2\}]\simeq P_3$, $I-N(u)$ is a stable set and hence $\chi(I-N(u))\le 1$.  Similarly, $\chi((I\cap N(u))-N(u'))\le 1$. Trivially, $(I\cap N(u)\cap N(u'))\cup \{u_2\}$ is stable set and hence the chromatic number is less than 2. Therefore, $\chi(G)\le \chi(I-N(u))+\chi((I\cap N(u))-N(u'))+\chi((I\cap N(u)\cap N(u'))\cup \{u_2\})\le 1+1+1=3.$

    $G[D_1\cup \{u_1\}\cup \{x\}\cup \{v_3\}]$ is a stable set. If $E[\{x\},D_1]\ne \emptyset$, then $G$ contains $co-domino2$. Accoding to definition, $E[\{u_1\},D_1\cup \{x\}]=\emptyset$. The rest is obvious.

    $G[D_2\cup D_4\cup (B_2-\{x\})]$ is a stable set. Recalling that $|C_3|\le 1$ and $C_1=C_2=\emptyset$, $B_2-\{x\}=C_4$. According to definition $D_2\cup D_4\cup (B_2-\{x\})$ is complete to $\{u_1,u_3\}$ and the rest is obvious.

    $G[D_3\cup \{u_3\}\cup \{v_1\}]$ is a stable set. According to definition, $E[D_3,u_3]=\emptyset$ and $G[E_3]$ is a stable set. The rest is obvious.

    Let $D_5'$ be $D_5\cup \{v_2\}$. $\chi(G)\le \chi(B_1\cup B_3-\{u_1,u_3\})+\chi(D_1\cup \{u_1\}\cup \{x\}\cup \{v_3\})+\chi(D_2\cup D_4\cup (B_2-\{x\}))+\chi(D_3\cup \{u_3\}\cup \{v_1\})\le 4+1+1+1=7.$

\begin{spacing}{2.0}
\end{spacing}

     Suppose $E[D_1,D_3]=\emptyset$.

     $G[D_1\cup D_3\cup \{x\}]$ is a stable set. If $E[\{x\},D_1\cup D_3]\ne \emptyset$, then $G$ contains $co-domino2$. Since $E[D_1,D_3]=\emptyset$, $G[D_1\cup D_3\cup \{x\}]$ .

     $\chi(G)\le \chi(D_1\cup D_3\cup \{x\})+\chi(D_2\cup D_4\cup (B_2-\{x\}))+\chi(B_1\cup \{v_2,v_3\})+\chi(B_3\cup \{v_1\})+\chi(D_5)\le 1+1+2+2+1\le 7$.

\begin{spacing}{2.0}
\end{spacing}

   Combining $co-domino2$-free, $do-domino3$-free and our supposition that evrey vertex in $A_2\cap N(v_2)$ has two neighbors in $\{u_1,u_2,u_3\}$, if $u\in N(v_1)\cap N(v_2)$, then $u$ is complete to $\{u_1,u_3\}$. Similarly, if $u\in N(v_2)\cap N(v_3)$ ,then $u$ is complete to $\{u_1,u_3\}$. Therefore, $G[(N(v_1)\cap N(v_2))\cup (N(v_2)\cap N(v_3))]$ is a stable set. 

   Furthermore,$\chi(N(v_2)-\{v_1,v_3\})\le 2$. Noticing that $N(v_2)=(N(v_1)\cap N(v_2))\cup (N(v_3)\cap N(v_2))\cup C_3\cup C_4$, the proof is easy.

   $\chi(G)\le \chi(N(v_2)-\{v_1,v_3\})+\chi(B_1\cup \{v_2,v_3\})+\chi(B_3\cup \{v_1\})+\chi(N(v_1)\cap N(v_3))\le 2+2+2+1=7.$

   Therefore, if $C_3\ne \emptyset$, $\chi(G)\le 7.$

\begin{spacing}{2.0}

\end{spacing}

    Suppose $C_1=C_2=C_3=\emptyset$. We continue to use notations $D_1,D_2,...,D_5$. Similar to case2, $\chi(G)\le 7$.(The proof does not depend on the existence on \{x\}.)\qed
\begin{spacing}{2.0}

\end{spacing}
 
Let $C = v_1v_2v_3u_3u_2u_1$ be a $6$-hole and $\{u\}$ be a vertex outside the $6$-hole. We use {\bf $X_1$} to denote a graph obtained from $\{u\}\cup C$ by connecting edges $v_1v_3,uv_1$ and $uu_2$. 

Let $C = v_1v_2v_3u_3u_2u_1$ be a $6$-hole and $\{u\}$ be a vertex outside the $6$-hole. We use {\bf $X_2$} to denote a graph obtained from $\{u\}\cup C$ by connecting edges $v_1v_3,v_3u_1$,$u_1u_3,v_2u$ and $uu_2$. 

The following figures are $X_1$ and $X_2$( from left to right).
 
\begin{figure}[H]
\centering
\includegraphics[width=0.5\textwidth]{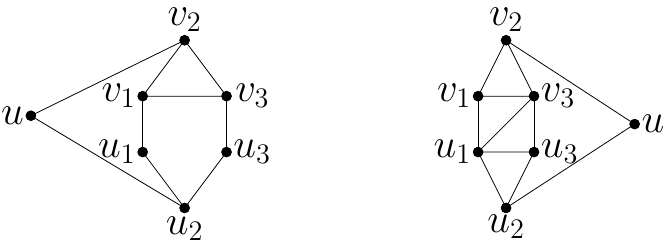}
\end{figure} 

\begin{lem}

Suppose $G$ is $co-domino$-free and contains $X_1$ or $X_2$, then $\chi(G)\le 7$.
\end{lem}

\pf

	{\bf case1:}G contains $X_1$.

     $G[B_1-A_0-\{u_1\}]$ is anticomplete to $\{u_1,u_2\}$ and complete to $\{u,u_3\}$.Suppose $y\in B_1-A_0$ and $y\sim u_1$(or $u_2$) only, then $G[\{y,u_1,u_2\}\cup \{v_3,v_1\}]\simeq P_3\cup P_2$. Therefore, $y$ must complete to both $u_2$ and $u_1$ and hence we get $G[\{y,u_1,u_2\}\cup \{v_1,v_3\}]\simeq K_3\cup P_2$,which contradics our assumption.   $G[B_2]$ is complete to $\{u_3\}$(or $u$) or $G[(B_2-A_0)\cup \{v_1,v_2\}\cup \{u_3,u_2\}]$ will induce an $P_3\cup P_2$($G[(B_2-A_0)\cup \{v_1,v_2\}\cup \{u,u_2\}]$ will induce a $P_3\cup P_2$).
  
    Symmetrically , $G[B_3-\{u_3\}]$ is anticomplete to $\{u_2,u_3\}$ and complete to $\{u,u_1\}$ and $G[B_2]-\{u\}$ is anticomplete to $\{u_1,u_3\}$ and anticomplete to $\{u,u_2\}$.
        
    Either $G[B_1-A_0]$ or $G[B_3]$ is edge-free. According to the above paragraph, $G[(B_1-A_0)\cup B_3]$ is complete to $u$ and hence it is $K_3$-free. Suppose there is an edge in $G[B_1-A_0]$($y_1,y_2$) and an edge in $G[B_3]$($y_3,y_4$). 
   
 We prove $G[v_4,y_1,y_2,y_3,y_4,v_6]\simeq co-domino$ to obtain a contradicion. It is easy to see that $y_3$ should have at least one neighbor in $\{y_2,y_1\}$. So does $y_4$. Consequently, $G[\{y_1,y_2,y_3,y_4\}]\simeq C_4$, otherwise $G[\{x_1,x_2,y_1,y_2\}]$ must contain $K_3$.  

   Symmetrically, either $G[B_1-A_0]$ or $G[B_2]$ is edge-free and either $G[B_2]$ or $G[B_3]$ is edge-free.

   $\chi(B_1\cup B_2\cup B_3)\le 4$. Recalling that $G[B_i]$ is $P_3$-free. Suppose $G[B_1-A_0]$ has edge, then $\chi(B_2)\le 1,\chi(B_3)\le 1$. If $G[B_2]$ or $G[B_3]$ contains edge, then we can prove $\chi(B_1\cup B_2\cup B_3)\le 4$. If $G[B_1-A_0]$,$G[B_2]$ and $G[B_3]$ are all edge-free, then it is trvial that $\chi(B_1\cup B_2\cup B_3)\le 4$.

    Therefore, $\chi(G)\le \chi(B_1\cup B_2\cup B_3)+\chi(A_2)\le 4+3=7.$

{\bf case2:} $G$ is $X_1$-free and contains $X_2$.

   Accoring to section 2.1, we partition $V(G)$ around $G[\{v_1,v_2,v_3\}]$.

   Every vertex in $B_2-\{u\}$ is adjacent  to $u_3$ and not adjacent  to $u_2$. If there is $x\in B_2-\{u\}$ such that $x$ is not adjacent to $u_3$, then $x$ is adjacent to $u_2$. Otherwise $G[\{v_1,v_2,x\}\cup \{u_2,u_3\}]\simeq P_3\cup P_2$ or $K_3\cup P_2$. Therefore, we can simply suppose there is $x\in B_2-\{u\}$ such that $x$ is adjacent to $u_2$. However, $G[\{u,u_2,x\}\cup \{v_1,v_3\}]\simeq P_3\cup P_2$ or $K_3\cup P_2$(depending on the adjacency of $u$ and $x$). Therefore, such $x$ does not exist.

   Every vertex in $B_1-A_0$ is adjacent to $\{u_3\}$. 

   Every vertex in $B_3$ is anticomplete to $\{u_2,u_3\}$.

   Partition $B_1\cup B_2\cup B_3-A_0$ into the folloing six sets.
\begin{itemize}
  \item $D_1:=\{x\in (B_2-\{u\})|x \text{ is adjacent to } u_1 \text{ and }u_3\}$.
  \item $D_2:=\{x\in (B_2-\{u\})|x \text{ is adjacent to } u_3 \text{ only}\}$.
  \item $D_3:=\{x\in (B_1-A_0)|  x \text{ is adjacent to } u_3 \text{ only}\}$.
  \item $D_4:=\{x\in (B_1-A_0)|  x \text{ is adjacent to } u_1 \text{ and }u_3\}.$
  \item $D_5:=\{x\in B_3|x \text{ is adjacent to }u_1\}.$
  \item $D_6:=\{x\in B_3|x \text{ has no neighbor in } \{u_1,u_2,u_3\}\}.$
\end{itemize}

  Since $D_1\cup D_4$ is compete to $\{u_1,u_3\}$ and $u_1\sim u_3$, $G[D_1\cup D_4]$ has no edge and hence $\chi(D_1\cup D_4)\le 1$. 

  Since $D_2$ and $D_3$ anticomplete to $\{v_3,u_1,u_2\}$, $G[D_2\cup D_3]$ is $P_2$-free and hence $\chi(D_2\cup D_3)\le 1$.

  $E[\{u\},A_0\cup D_6 ]=\emptyset$. Since $\{u\}\cup D_6$ is anticomplete to $\{v_1,u_1,u_3\}$ and $G[\{v_1,u_1,u_3\}]\simeq P_3$, $G[\{u\}\cup D_6]$ has no edge. Since $\{u\}\cup A_0$ is anticomplete $v_1,v_3,u_3$ and $G[\{v_1,v_3,u_3\}]\simeq P_3$(It is trivial that $u_3$ is anticomplete to $A_0$), $G[\{u\}\cup A_0]$ has no edge.

  Since $D_6\cup A_0$ is anticomplete to $\{v_2,v_1,u_1\}$ and $G[\{v_2,v_1,u_1\}]\simeq P_3$, $G[A_0\cup D_6]$ has no edge and hence $G[\{u\}\cup D_6\cup A_0]$ has no edge. 

  $\chi(B_1\cup B_2\cup B_3)\le \chi(D_1\cup D_4)+\chi(D_2\cup D_3)+\chi(A_0\cup \{u\}\cup D_6)+\chi(D_5)\le 4$. Therefore, $\chi(G)\le \chi(B_1\cup B_2\cup B_3)+\chi(A_2)\le 7.$\qed
\begin{spacing}{2.0}

\end{spacing}
Suppose all graphs are $(co-domino,K_3\cup P_2,X_1,X_2)$-free.

 Let $C = v_1v_2v_3v_4v_5$ be a $5$-hole. We use {\bf $co-twin-C_5$} to denote a graph obtained from $C$ by adding a vertex $\{v_6\}$ which is only adjacent $\{v_3,v_4,v_5\}$.
\begin{figure}[H]
\centering
\includegraphics[width=0.20\textwidth]{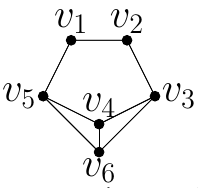}
\end{figure}

    Let $C = v_1v_2v_3u_3u_2u_1$ be a $6$-hole and $\{u\}$ be a vertex outside the $6$-hole. We use {\bf $\mathcal{Y}$} to denote a family of graphs obtained from $\{u\}\cup C$ by connecting edge $uv_1,uv_2$ and $uu_2$ and it does not matter whether $uu_1$ is connected or not.

\begin{figure}[H]
\includegraphics[width=0.20\textwidth]{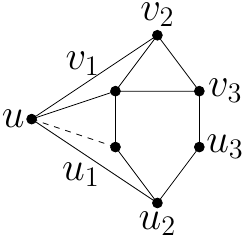}
\end{figure}
   Define $M(v_1,v_2):=G-N(v_1)-N(v_2)-\{v_1,v_2\}$. In this lemma, $N(v_1)$ does not include $\{v_2\}$ and $N(v_2)$ does not include $\{v_1\}$. We introduce a partition which is useful in the following lemma.
\begin{obs}{\label{o4}}
   If $G$ contains $co-twin-C_5$. $V(G)$ can be partitioned into $\{v_1,v_2\},N(v_2)\cup N(v_1)$  and $M(v_1,v_2)$. 
\end{obs}
                                          
\begin{lem}
   If $G$ contains $co-twin-C_5$, then $\chi(G)\le 7$.
\end{lem}
\pf

{\bf case1:} G contains an element of $\mathcal{Y}$.

   $G[B_1-A_0]$ is complete to  $\{u_3\}$ and anticomplete to $\{u_1,u_2\}$. Otherwise, suppose there is $x\in G[B_1-A_0]$ which is not adjacent to $\{u_3\}$ or has neighbor in  $\{u_1,u_2\}$. If $x$ has neighbor  in $\{u_1,u_2\}$, then $G[\{x,u_1,u_2\}\cup \{v_2,v_3\}]$ induces $P_3\cup P_2$ or $K_3\cup P_2$. Therefore, $x$ is not adjacent to $u_3$ and hence $G[\{v_2,x\}]\simeq P_2$ is anticomplete to $G[\{u_1,u_2,u_3\}]\simeq P_3$. Therefore, such $x$ does not exist.

   Similarly, $G[B_3]$ is complete to $\{u_1\}$ and anticomplete to $\{u_3,u_2\}$.

   Either $G[B_1-A_0]$ or $G[B_3]$ has no edge.  Otherwise, suppose there is an edge($y_1,y_2$) in $G[B_1-A_0]$ and an edge($y_3,y_4$) in $G[B_3]$.  If $G[\{y_1,y_2,y_3,y_4\}]$ induces $diamond$ or $C_4$, then $G[\{u_1,u_2,u_3,y_1,y_2,y_3,y_4\}]$ induces $X_1$ or $X_2$. Therefore, $G[\{y_1,y_2,y_3,y_4\}]$ is not isomorphic to $diamond$, $C_4$ and $K_4$. However, such condition requires $G[\{y_1,y_2,y_3,y_4\}]$ to be isomorphic to $P_4$ or $2K_2$. Therefore, there must be a vertex has degree one in $G[\{y_1,y_2,y_3,y_4\}]$. Suppose $d(y_1)=1$ in $G[\{y_1,y_2,y_3,y_4\}$. Then $y_1$ is anticomplete to $\{y_3,y_4\}$ and hence $G[\{u_2,u_3,y_1\}\cup \{y_3,y_4\}]\simeq P_3\cup P_2$. Therefore, either $y_1,y_2$ does not exist or $\{y_3,y_4\}$ does not exists.

  $G[B_2]$ is complete to $\{u_3\}$. Otherwise, suppose there exists $z\in G[B_2]$ such that z is not adjacent to $\{u_3\}$. $z$ is adjacent to $u_2$, otherwise $G[\{v_1,v_2,z\}\cup \{u_2,u_3\}]\simeq P_3\cup P_2$. Since $G[\{x,u_1,u_2,u_3,v_1,v_2,v_3\}]\simeq X_1$, $z$ is adjacent to $u_1$. However, if $z$ is adjacent to $u_1$, then $G[\{z,u_1,u_2,v_1,v_2,v_3\}]\simeq co-domino$. Therefore, such $z$ does not exist.

  $G[B_2]$ contains at most one edge. Suppose there are two edges in $G[B_3]$($\{y_1,y_2\}$ and $\{y_3,y_4\}$ ). Because $v_2$ is complete to $\{u,y_1,y_2,y_3,y_4\}$. $u$ has at most one neighbor in $\{y_1,y_2\}$ and $\{y_3,y_4\}$. Suppose $u\not\sim y_2,u\not\sim y_3$. Since both $y_2,y_3$ are not adjacent to $v_1$, $G[\{y_3,u_3,y_2\}\cup \{u,v_1\}]\simeq P_3\cup P_2$, which causes contradiction. Therefore, $G[B_2]$ has at most one edge.

  Suppose $B_1-A_0$ has no edge. Let $B_2=V_1, B_3=V_2,A_0=V_3$, then we can apply claim\ref{c1} to get $\chi(B_1\cup B_2\cup B_3)\le \chi(B_1-A_0)+\chi(B_2\cup B_3\cup A_0\le 4$ and hence $\chi(G)\le \chi(B_1\cup B_2\cup B_3)\le 4+3=7$.

{\bf case2:} G is $\mathcal{Y}$-free and contains a $co-tiwn-C_5$.

   According to observation \ref{o4}, $V(G)=\{v_1,v_2\}\cup(N(v_2)\cup N(v_1))\cup  M(v_1,v_2)$.

   Noticing that vertices in $N(v_1)\cup N(v_2)-\{v_3,v_5\}$ have neighbor in $\{v_4,v_6\}$, we can divide $N(v_1)\cup N(v_2)-\{v_3,v_5\}$ to make the structure more clearly.  We difine two disjoint sets as following:
\begin{itemize}
\item $D_1:=\{y|y \in N(v_1)\cup N(v_2)-\{v_3,v_5\},y$ is adjacent to only one vertex in $\{v_3,v_4,v_5,v_6\}\}$
\item $D_2:=\{y|y \in N(v_1)\cup N(v_2)-\{v_3,v_5\},y$ is adjacent to more than one vertices in $\{v_3,v_4,v_5,v_6\}\}$
\end{itemize}   

  $|D_1|=\emptyset$. Otherwise, suppose $y_1\in D_1$ and $y_1\sim v_1$. Suppose $y_1$ is adjacent to $\{v_4\}$.
Because if the only neighbor of $y_1$ belongs to $\{v_3,v_5\}$ , then $G[\{v_1,v_2,y_1,v_4,v_6\}]$ induces $K_3\cup P_2$ or $P_3\cup P_2$.
 
   If $y_1$ is not adjacent to $v_2$, then $G[\{v_1,y_1,v_2,v_3,v_4,v_5,v_6\}]$ is isomorphic to an element in  $\mathcal{Y}$(when $u_1$ is not adjacent to $u$. However, if $y_1$ is not adjacent to $v_2$, then $G[\{v_1,y_1,v_2,v_4,v_5,v_6\}]\simeq co-domino $. Therefore, such $y_1$ does not exist.

    Any vertex in $D_2$ is complete to an edge in $G[\{v_3,v_4,v_5,v_6\}]$. Otherwise, there is $z\in D_2$ which is not complete to any edge in $G[\{v_3,v_4,v_5,v_6\}]$. $z$ must be adjacent to $\{v_3,v_5\}$, which leads to the contradiction that $G[\{v_1,v_2,z,v_4,v_6\}]$ induces $K_3\cup P_2$ or $P_2\cup P_3$. 

   $\{v_1,v_2\}\cup M(v_1,v_2)$ can be colored  with no more than $2$ colors.

    Recalling that points in $D_2$ are at least adjacent to one edge in $G[\{v_3,v_4,v_5,v_6\}]$, it is trivial that $\chi(D_2)\le 5$. In summary, $\chi(G)\le \chi(\{v_1,v_2\}\cup M(v_1,v_2))+\chi(D_2\cup \{v_3,v_4,v_5,v_6\})\le 2+5\le7$.        \qed  

\begin{spacing}{2.0}
\end{spacing}

Suppose all graphs are $(co-domino,K_3\cup P_2,X_1,X_2,co-twin-C_5)$-free. 
  
   Let $C = u_1v_1v_2v_3u_3u_2$ be a $6$-hole. We use {\bf $\chi37$} to denote a graph obtained from $C$ by connecting edge $v1v3$.
\begin{figure}[H]
\centering
\includegraphics[width=0.40\textwidth]{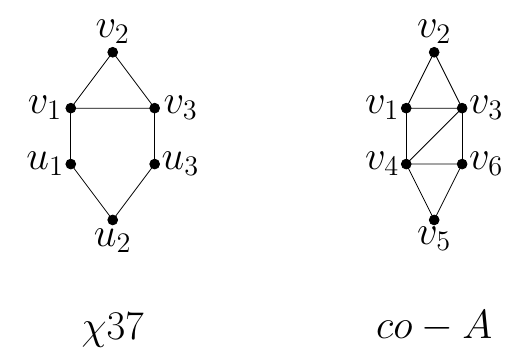}

\end{figure}

\begin{lem}{\label{l2}}
    If $G$ contains $\chi37$ or $co$-$A$, then $\chi(G)\le 7$.
\end{lem}
\pf
 
{\bf case1:} $G$ contains $\chi37$.

    We prove $G[B_2]$ is edge-free. We obtain contradiction by proving $G[\{v_1,v_3,u_1,u_3,x_1,x_2\}]\simeq co-twin-C_5$. In other words, contradiction comes out if $x_1,x_2$  are anticomplete to $\{u_2\}$ and complete to $\{u_1,u_3\}$. Suppose there is an edge in $G[B_2]$ and we call them $G[x_1,x_2]$.  If $x_1$(or $x_2$), $G[\{v_1,v_3\}\cup \{x_1,x_2,u_2\}]$ would induce $K_3\cup P_2$ or $P_3\cup P_2$. Furthermore, we can  prove $\{x_1,x_2\}$ is complete to $\{u_1,u_3\}$ without much difficulty and hence we finish our proof of proposion that $G[B_2]$ is edge-free.

    Furthermore, either $G[B_3]$ or $G[B_1-A_0]$ has no edge. Otherwise, suppose there is $\{y_1,y_2\}\in E(B_1-A_0)$, $\{y_3,y_4\}\in E(B_3)$. Since $B_1-A_0$ is complete  to $\{u_3\}$ and anticomplete to $\{u_1,u_2\}$ and $B_3$ is complete to $\{u_1\}$ and anticomplete to $\{u_2,u_3\}$, $G[y_1,y_3,u_1,u_3,u_2,y_2]$ induces $co-domino$ or $co-twin-C_5$. 
   
    According to claim\ref{c1}, $\chi(B_1\cup B_3)\le 3$. Therefore, $\chi(G)\le \chi(B_2)+\chi(B_1 \cup B_3)+\chi(A_2)\le 7$.
\begin{spacing}{2.0}
\end{spacing}

    {\bf case2:} $G$ is $\chi37$-free and contains a $co$-$A$.

    $G[B_2]$ is complete to $\{v_4,v_6\}$ and hence $G[B_2]$ is edge-free. Otherwise, there is $x\in B_2$ such that  $x$ has only one neighbor in $\{v_4,v_6\}$. If $x\sim v_4$, then $x\sim v_5$. Otherwise, $G[\{v_1,v_2,x\}\cup \{v_5,v_6\}]$ induces $P_3\cup P_2$. However, $G[\{x,v_2,v_3,v_6,v_5,v_1\}]$ induces $\chi37$. Therefore, $x\sim v_6$. Now $v\not\sim v_5$, otherwise $G[\{v_1,v_2,v_3,x,v_6,v_5\}]$ induces $co-domino$. However, $G[\{v_1,v_2,x,v_6,v_4,v_5\}]$ induces $co-A$. Therefore, such $x$ does not exist.

   every edge in $G[B_1-A_0]$ includes one vertex which is complete to $\{v_4,v_6\}$. Otherwise, there is an edge in $G[B_1-A_0]$ whose vertices has at most one neighbor in $\{v_4,v_6\}$. If one of them($x$) is adjacent to $\{v_4\}$ rather than $\{v_6\}$,then  $G[\{v_2,v_1,x\}\cup \{v_5,v_6\}]$ induces $P_3\cup P_2$. Therefore, none of these vertices is adjacent to $\{v_4\}$ and hence $G[\{v_3,v_4,v_5\}\cup B_1-A_0]$ induces $P_3\cup P_2$. We already obtain a contradiction.

   According to the above discussion, $G[(B_1-A_0)\cup B_2]$ can be partitioned into two stable sets: the vertices which are complete to $\{v_4,v_6\}$ and the rest of vertices in $B_1-A_0$. Therefore, $\chi((B_1-A_0)\cup B_2)\le 2$.

   Therefore, $\chi(G)\le \chi(B_2\cup (B_1-A_0))+\chi(B_3\cup A_0)+\chi(A_2)\le 2+2+3=7.$\qed

\begin{spacing}{2.5}

\end{spacing}

Combining lemma\ref{l1} and lemma\ref{l2}, we obtain theorem\ref{t2} straightly. 
\begin{thmm}{1.2}
      If $G$ contains $co-domino$ or $co-A$, then $\chi(G)\le 7$.
\end{thmm}
  Suppose $G$ is $(K_3\cup P_2,co-domino,co-A)$-free with $\overline{G[D_1]}$ is $P_3$-free.

\begin{claim}{\label{c8}}
  If $v_1,v_2\in D_1$ and $v_1\not\sim v_2$, then either $G[N(v_1)-N(v_2)]$ or $G[N(v_2)-N(v_1)]$ has no edge.
\end{claim}
\pf
     
      Suppose both sets have edges. The vertex sets in $G[N(v_1)-N(v_2)]$ is $\{u_1,u_2\}$ and that of edge in $G[N(v_2)-N(v_1)]$ is $\{u_3,u_4\}$. $G[\{v_1,v_2,u_1,u_2,u_3,u_4\}]$ must induce $co-domino$ or $co-A$ or $P_3\cup P_2$. Therefore, either $G[N(v_1)-N(v_2)]$ or $G[N(v_2)-N(v_1)]$ is edge-free.
 \qed

\begin{thmm}{(1.3)}
       If $G[D_2]$  is not a clique, then $\chi(G)\le $7.
\end{thmm}

\pf   

     Since $G[D_2]$ is not a clique, $G[D_2]$ has two  nonadjacent vertices $\{v_1,v_2\}$. According to claim\ref{c8}, either $G[N(v_1)-N(v_2)]$ or $G[N(v_2)-N(v_1)]$ has no edge. Without loss of generality, we suppose $G[N(v_1)-N(v_2)]$ is edge-free and hence $\chi(N(v_1)\cup N(v_2))\le 1$.
\begin{spacing}{2.0}
\end{spacing}

     If $\omega(N(v_1)\cap N(v_2))\le 1$, then $\chi(N(v_1)-N(v_2))\le 1$ or $\chi(N(v_2)-N(v_1))\le 1$. We suppose $\chi(N(v_1)-N(v_2))  \le 1$, then $\chi(v_1)\le 2$. If $G[N(v_1)]$ has edge($\{v_2,v_3\}$), then accoring to section 2.1, $G[\{v_1,v_2,v_3\}]$ can be divided into $B_1\cup B_2\cup B_3\cup (N(v_1)\cap N(v_2))\cup (N(v_1)\cap N(v_3)) \cup (N(v_2)\cap N(v_3))$. 

    $\chi(N(v_1))\le \chi(N(v_1)\cap N(v_2))+\chi(N(v_1)-N(v_2))\le 2$. The vertex sets left to be colored are:$A_0,N(v_2)\cap N(v_3),B_2$ and $B_3$. Since $\chi(B_2\cup A_0)\le 2$, $\chi(B_3)\le 2$ and $\chi(N(v_2)\cap N(v_3))\le 1$. Therefore,  $\chi(G)\le 7.$ 
\begin{spacing}{2.0}
  
\end{spacing}
     Suppose $\omega(N(v_1)\cap N(v_2))=2$.
      
      We partition $V(G)$ into three parts: $\{v_1,v_2\}$, $N(v_1)\cup N(v_2)$ and $G-\{v_1,v_2\}-N(v_1)-N(v_2)$. Define $A:=G-\{v_1,v_2\}-N(v_1)-N(v_2)$.

      $\chi(A)\le 3$. Suppose there is an edge($x_1x_2$) in $G[A]$. Then $A$ can be divided into $A-N(x_1)$ and  $(A-N(x_2)) \cap N(x_1)$ and $A\cap N(x_1) \cap N(x_2)$. Since $A-N(x_1)$ is anticomplete to $\{v_1,x_1,v_2\}$ and $G[\{v_1,x_1,v_2\}]\simeq P_3$, $G[A-N(x_1)]$ must be $P_2$-free and hence $\chi(A-N(x_1))\le 1$. Similarly, $\chi((A-N(x_2)) \cap N(x_1))\le 1$. Because $x_1\sim x_2$ and $\omega(G)\le 3$, $\chi(A\cap N(x_1) \cap N(x_2))\le 1$.  Therefore ,$\chi(A)\le 3$. 

     $\chi(N(v_2))\le 3$. According to definition of $D_2$, there is a $K_3(\{u_1,u_2,u_3\})$ induced in  $G[G-N(v_2)]$. Define $E_{i}:=\{x\in N(v_2)|x \text{ only adjacent to }\{u_i\}\}$ and $E_{i,j}:=\{x\in N(v_2)| x \text{ only adjacent to } \{u_i\} \text{ and } \{u_j\}\}$. Since $G$ is $co$-$A$-free, $E_i$ is anticomplete to $E_{i,i+1}$(mod$3$). Obviously, $|E_i|\le 1$ for $i=1,2,3$ and $G[E_{i,j}]$ is edge-free. Therefore $\chi(N(v_2))\le \sum_{i=1}^{3}\chi(E_{i}\cup E_{i,i+1})\le 3$. 

       Because $\chi(N(v_2))\le 3$, $\chi(N(v_2)\cup N(v_1))\le \chi(N(v_1)-N(v_2))+\chi(N(v_2))\le 1+3=4$.  Consequently, $\chi(G)\le \chi(N(v_1)\cup N(v_2))+\chi(\{v_1,v_2\}\cup A)\le 4+3$=7.       \qed                      

Therefore, combining theorem\ref{t1}, theorem\ref{t2} and theorem\ref{t3}, we obtain our major theorem.

\begin{thmm}
   If $G$ is $(P_3\cup P_2, K_4)$-free, then $\chi(G)\le 7$.
\end{thmm}

  Finally, we prove a simple theorem using the above theorem.

\begin{thmm}{(1.4)}
   If $G$ is $(4K_1,\overline{P_3\cup P_2})$-free with order $n$ and clique number $\omega$, then $n\le 7\omega$ and $\chi(G)\le 4\omega$.
\end{thmm}
\pf

   Suppose $G$ is $(4K_1,\overline{P_3\cup P_2})$-free.

   If $\overline{G}$ is connected, then we apply the above theorem to obtain that $\chi(\overline{G})\le 7$. Therefore, $V(\overline{G})$ can be partitioned into $7$ stable sets and hence $V(G)$ can be partitioned into $7$ cliques. Therefore, $n\le 7\omega$.

   If $\overline{G}$ is not connnected, then there is $G_1$ such that $G=G_1+G_2$ and hence $|V(G)|=|V(G_1)|+|V(G_2)|$. Such decomposition will last until each $G_i$ satidfied $\overline{G_i}$ is connected. Suppose $G=G_1+...G_k$. Since $V|G_i|\le 7\omega(G_i)$ and $\omega=\sum_{j=1}^k\omega(G_i)$, $|V(G)|=n\le 7\omega$.

   It is trivial that the complement of bipartite graph is perfect graph. Therefore $V(G)$ can be partitioned into $3$ perfect graphs and one clique. Therefore, $\chi(G)\le 4\omega$.\qed

\section{acknowledge}

\quad I try to go further on $(P_3\cup P_2,K_4)$-free graphs. However, I fail and 7 is the best upper bound I obtain.  I know this is not a well-written paper. If you are confused about some details in the proof, you can email me at 1345770246@qq.com or SA23001024@mail.ustc.edu.cn.
    
\end{document}